\documentclass[a4paper,twosided]{article}

\usepackage[a4paper,
            left=4cc,
            right=4cc]{geometry}

\usepackage{times,cite}
\usepackage[T1]{fontenc}
\usepackage[normalem]{ulem}
\usepackage{hyperref}
\usepackage[utf8]{inputenc}
\usepackage{graphicx}
\usepackage{amsmath}
\usepackage{amsfonts}
\usepackage{authblk}
\usepackage{xcolor}
\providecommand{\WileyBibTextsc}{}
\let\textsc\WileyBibTextsc
\providecommand{\othercit}{}
\providecommand{\jr}[1]{#1}

\begin{document}

\title{A brief review of Reduced Order Models using intrusive and non-intrusive techniques}

\author[1]{Guglielmo~Padula\footnote{gpadula@sissa.it}}
\author[1]{Michele~Girfoglio\footnote{mgirfogl@sissa.it}}
\author[1]{Gianluigi~Rozza\footnote{grozza@sissa.it}}

\affil[1]{Mathematics Area, mathLab, SISSA, via Bonomea 265, I-34136 Trieste,
Italy}

\maketitle                  
\begin{abstract}
Reduced Order Models (ROMs) have gained a great attention by the scientific community in the last years thanks to their capabilities of significantly reducing the computational cost of the numerical simulations, which is a crucial objective in applications like real time control and shape optimization.  This contribution aims to provide a brief overview about such a topic.  We discuss both an intrusive framework based on a Galerkin projection technique and non-intrusive approaches, including Physics Informed Neural Networks (PINN), purely Data-Driven Neural Networks (DDNN), Radial Basis Functions (RBF), Dynamic Mode Decomposition (DMD) and Gaussian Process Regression (GPR). We also briefly mention geometrical parametrization and dimensionality reduction methods like Active Subspaces (AS). Then we present some results related to  academic test cases as well as a preliminary investigation related to an industrial application.
\end{abstract}

\section{Introduction}
The big power of modern computer systems has of course enabled the development and the employment of complex numerical methods. However, in certain cases, including real time control of systems like digital twins \cite{democont, digtwin} or optimization problems \cite{shapeopt, shapeopt2}, it is still very hard to obtain an high computational efficiency.  \\
There are mainly two ways to solve this problem. The first one is to use High Performance Computing (HPC) clusters \cite{hpc} but they are immobile, costly to build and have high energy requirements. The second is to introduce Reduced Order Models (ROMs) enabling fast computation while introducing a slight loss of accuracy.\\
The development of a ROM is typically split in two stages:
\begin{itemize}
    \item The offline stage, in which few selected values of the parameters (physical or geometrical) involved in the problem, high-fidelity simulations are computed, typically by using an HPC facility. 
    Such simulations, obtained by solving the discretized governing equations system, hereinafter referred to as Full Order Model (FOM), are computationally expensive, but they are performed only once. At this point one makes use of some techniques, among which Proper Orthogonal Decomposition (POD) \cite{hesthaven,rozza}, to extract the essential information from the high-fidelity database and to assemble the surrogate model.
    \item The online stage, in which, taking advantage of the surrogate model, the solution is efficiently predicted for new values of the parameters. The computational speedup with respect to the FOM could be significantly high: for instance, for some fluid dynamics problems, a value of $1e5$ is achievable (see Sec. \hyperref[sec:3]{3}). 
\end{itemize}
ROMs are divided in two big families\cite{rozza}: intrusive methods, in which one manipulates directly the governing equations, and non-intrusive methods, in which only the simulation data are considered. It should be noted that also for non-intrusive methods one could make use of information coming from the physics
of the system at hand, as it happens for PINN (see Sec. \ref{sec:2}). \\
The paper is organized as follows. In Sec.  \hyperref[sec:2]{2}, an intrusive technique, the Galerkin projection method, and several data-driven techniques,  Physics Informed Neural Networks (PINN), Radial Basis Functions (RBF), purely Data-Driven Neural Network (DDNN), Dynamic Mode Decomposition (DMD) and Gaussian Process Regression (GPR) are briefly discussed. The geometrical parametrization is also introduced and the Active Subspace (AS) method is presented as a possible solution strategy. In Sec. \hyperref[sec:3]{3} we investigate some benchmark test cases and in Sec. \hyperref[sec:4]{4} we present an industrial application. 
Finally in Sec. \hyperref[sec:5]{5} we report some concluding remarks.

\section{Reduced Order Models}
\label{sec:2}
We refer to ROMs to accelerate the numerical solution of parameterized partial differential equations (PDEs) in the space time domain, whose general form is the following: 
\begin{equation}
\begin{cases}
\mathcal{L}(u,\mu, t, x, y, z)=0 \quad \text{with} & (x, y, z) \in \Omega(\mu), \quad t \in(0, T], \quad \mu\in \mathbb{P}, \\
\mathcal{B}(u,\mu, t, x, y, z)=0 \quad \text{with} & (x, y, z) \in \partial \Omega(\mu), \quad t \in(0, T], \quad \mu \in \mathbb{P}, \\
\mathcal{T}(u,\mu,x,y,z)=0 \quad \text{with} & (x, y, z) \in \Omega(\mu), \quad t=0, \quad \mu \in \mathbb{P}, 
\end{cases}
\label{eq:1}
\end{equation}
where  $\Omega(\mu)\subset \mathbb{R}^{3}$ is the space domain of the problem, $T$ is the final time of the simulation, $\mathbb{P}$ is the parameter space, $u\in \mathbb{V}$ and $\mathbb{V}\subseteq \{v | v: \mathbb{P}\times \mathbb{R}^{3}\rightarrow \mathbb{R}^{n}\}$ is an Hilbert space. Finally $\mathcal{L}:\mathbb{V}\times \mathbb{P}\times\mathbb{R}^{+}\times \mathbb{R}^{3}\rightarrow\mathbb{R}$, $\mathcal{B} :\mathbb{V}\times \mathbb{P}\times\mathbb{R}^{+}\times \mathbb{R}^{3}\rightarrow\mathbb{R}$, $\mathcal{T}:\mathbb{V}\times \mathbb{P}\times{\mathbb{R}}^{3}\rightarrow \mathbb{R}$ are continuous functions whose explicit formulation depends on the problem at hand.

\subsection{An intrusive method: the Galerkin Projection}
For sake of simplicity, we assume that:
\begin{itemize}
    \item the problem is steady, so $\mathcal{L}$, $\mathcal{B}$ and $\mathcal{T}$ do not depend on $t$.
    \item $\Omega$ does not depend on $\mu$.
    \item $\mathbb{V}=\{v\in \mathbb{P}\times \mathbb{R}^{3} | v(\mu,\cdot)\in H^{1}(\mathbb{R}^{3}) \text{ and } v(\mu,\cdot )_{\partial \Omega }=0 \text{ } \forall \mu \in \mathbb{P}\}$, so in this case $n=1.$
    \item The first equation of system \eqref{eq:1} can be rewritten as 
    \begin{equation}
    \sum_{i=1}^{m}\theta_{\mathcal{L}}^{i}(\mu)\mathcal{L}_{i}(u, x, y, z)=\sum_{i=1}^{m}\theta_{f}^{i}(\mu)f_{i}(x, y, z) \
    \quad \forall (x,y,z)\in \Omega, \mu \in \mathbb{P},
    \end{equation}
where $\theta_{\mathcal{L}}^{i}$ and $\theta_{f}^{i}$ are scalar continuous functions from $\mathbb{P}$ to $\mathbb{R}$, $f_i$ are continuous functions from $\mathbb{R}^{3}$ to $\mathbb{R}$ and $\mathcal{L}_i$ are continuous functions from $\mathbb{V}\times \mathbb{R}^{3}$ to $\mathbb{R}$.
    \item The second equation of system \eqref{eq:1} reads
    \begin{equation}
        u(\mu,x,y,z)=0 \quad \forall  (x,y,z)\in \partial\Omega, \mu \in \mathbb{P}.
    \end{equation}
\end{itemize}
By multiplying  by $v \in \mathbb{W}=\{v:\mathbb{R}^{3}\rightarrow \mathbb{R}| v_{\partial \Omega}=0\}$ and integrating over $\Omega$ we obtain the weak form of the problem:
\begin{equation}
    \sum_{i=1}^{m}\theta_{\mathcal{L}}^{i}(\mu)\mathcal{A}_{i}(u(\mu,\cdot),v)=\sum_{i=1}^{m}\theta_{f}^{i}(\mu)\mathcal{F}_{i}(v) \quad \forall v \in \mathbb{W},\mu \in \mathbb{P},
\label{eq:2}
\end{equation}
where 
\begin{equation}
\mathcal{A}_{i}(u(\mu,\cdot),v)=\int_{\Omega}u(\mu,x,y,z)v(x,y,z)dxdydz \quad \text{and} \quad \mathcal{F}_{i}(v)=\int_{\Omega}f_{i}(x,y,z)v(x,y,z)dxdydz.
\end{equation}
Based on our assumptions, $\mathcal{A}_{i}$ are bilinear with respect to $u,v$. Moreover, we also assume that they are symmetric and coercive, i.e. 
\begin{equation}
\mathcal{A}_{i}(u,v)>c_{i}||u||_{W}||v||_{W} \quad c_{i}>0 \quad  \forall u,v\in \mathbb{W}, \forall i = 1, 2, \dots, m.
\end{equation}
In this setting, it could be proven that a solution $u\in \mathbb{V}$ exists and is unique.
As the $\mathcal{A}_{i}$ are bilinear we can build the space $\mathbb{W}_{N}=\operatorname{span}\{v_{1},\dots,v_{N}\}.$ We can restrict to $u \in \{u\in \mathbb{V}|u(\mu,\cdot)\in \mathbb{W}_{N}, \text{  } \forall \mu \in \mathbb{P}$\}, so it follows that 
\begin{equation}
u=\sum_{j=1}^{N}\alpha_{j}(\mu)v_{j}, \label{eq:reduction}
\end{equation}
where $\alpha_{j}:\mathbb{P}\rightarrow \mathbb{R}.$ By substituting eq. \eqref{eq:reduction} in eq. \eqref{eq:2} we obtain the following system of $N$ equations 
\begin{equation}
\sum_{i=1}^{m}\sum_{j=1}^{N}\alpha_{j}(\mu)\theta_{\mathcal{L}}^{i}(\mu)\mathcal{A}_{i}(v_{j},v_{k})=\sum_{i=1}^{m}\theta_{f}^{i}(\mu)\mathcal{F}_{i}(v_{k})\quad k=1,\dots, N.
\end{equation}
Note that once a basis of $\mathbb{W}_{N}$ has been chosen, the computation of the $\mathcal{A}_{i}$ and $\mathcal{F}_{i}$ (which can be very costly) is done only once. 
We assume that it does exist a reduced space $\mathbb{W}_{rb}$ of dimension $N_{rb}<<N$  such that the error
\begin{equation}
\varepsilon = ||u(\mu,\cdot)-u_{rb}(\mu,\cdot)||_{\mathbb{W}}
\end{equation}
is small. The Galerkin projection method provides a posterior error bound for coercive problems \cite{hesthaven}. The methods used for computing the reduced basis are mainly two: the Proper Orthogonal Decomposition (POD) and the Greedy Algorithm \cite{hesthaven}. 
A disadvantage of the Galerkin projection approach is that, for complex problems, like Navier-Stokes equations\cite{prusak,girfoglio,sheidani,mola,khamlich}, there are no error bounds and the error decreases slowly. 

\subsection{Non-intrusive methods}
We briefly introduce the following techniques: Physics Informed Neural Network (PINN) \cite{democont,rozza}, purely Data-Driven Neural Networks (DDNN) \cite{rozza,padula,coscia,prusak,salavatidezfouli,meneghetti,hinze, gonnella,siena,regazzoni,papapicco,ivagnes,sienalibro}, Radial Basis Functions (RBF) \cite{rbf,padula,shapeopt3,coscia,hajisharifi, balzotti,ivagnes,sienalibro}, Dynamic Mode Decomposition (DMD) \cite{rozza,hajisharifi,hinze,andreuzzi} and Gaussian Process Regression (GPR) \cite{rozza,padula,shapeopt3,ivagnes}.

\subsubsection{Physics Informed Neural Networks}
A fully connected neural network with lengths $l_{0},\text{ } l_{1}, \dots, l_{L}$ is a parametrized function $h: x \in \mathbb{R}^{l_{0}}\rightarrow y \in \mathbb{R}^{{l_{L}}}$ which has the following form 
\begin{equation}
\begin{cases}
x_{0}=x,\\
x_{i}=f(A_{i}x_{i-1}+b_{i}) & \text{ with } i=1\dots L-1 \\
x_{L}=g(A_{L}x_{l-1}+b_{L}),\\
y=x_{L},
\end{cases}
\end{equation}
where 
$A_{i}\in \mathbb{R}^{l_{i-1}\times l_{i}},\text{  }$ $ b_{i}\in \mathbb{R}^{l_{i}}$ are learnable parameters and $f$ is a fixed non polynomial function called activation function. The function $g$ is typically the identity for regression problems and the \emph{softmax} function for classification problems \cite{rozza}. Neural networks are useful as they provide an approximation $u_\theta$ of any continuous function $u$ on compact sets:
\begin{equation}
u_{\theta}:\mathbb{R}^{dim(\mathbb{P})+4}\rightarrow \mathbb{R}^{n}.
\end{equation}
A Physics Informed Neural Network (PINN) \cite{democont,rozza} is a neural network which is trained using the following loss function:
\begin{eqnarray}
    \sum_{\mu \in \mathbb{P}}\sum_{(x,y,z)\in \Omega(\mu)}\sum_{t\in (0,T]} ||\mathcal{L}(u_{\theta}(\mu,t,x,y,z),t,x,y,z)||+
    \sum_{\mu \in \mathbb{P}}\sum_{(x,y,z)\in \partial\Omega(\mu)}\sum_{t\in (0,T]} ||\mathcal{B}(u_{\theta}(\mu,t,x,y,z),t,x,y,z)||+ \nonumber\\\sum_{\mu \in \mathbb{P}}\sum_{(x,y,z)\in \Omega(\mu)}||\mathcal{T}(u_{\theta}(\mu,0,x,y,z),\mu,x,y,z)||.
    \end{eqnarray}
The sum above is computed using Monte Carlo sampling\cite{rozza,democont}. In the offline phase, the neural network is trained using the backpropagation algorithm \cite{democont, rozza}. 
\subsubsection{Purely Data-Driven Neural Networks}
After that the full order database has been computed, we can approximate $u$ as a function  $u_{\theta}:\mathbb{R}^{dim(\mathbb{P})}\rightarrow \mathbb{R}^{M\cdot n \cdot N_{T}}$, where $N_{T}$ is the number of time instants, $n$ is the dimension of the output, and $\theta$ are the parameters of the neural network. To train it we adopt a simple $L_2$ loss
\begin{equation}
\min_{\theta\in {\mathbb{R}^{dim(\theta)}}}||u-u_{\theta}||_{2}^{2}.
\end{equation}
In this case when fed with a new parameter in the online phase, the neural network will output the entire solution on the discretized domain. This approach is typically faster than PINN one but does not support the computation of the solution in points different from the ones of the original discretization.

\subsubsection{Radial Basis Function}
Radial Basis Function (RBF) is an interpolation method that has the following form
\begin{equation}  
f(x)=\sum_{i=1}^n w_i \psi\left(\left\|x-x_i\right\|\right)+\sum_{j=1}^{m}c_{j}p_{j}(x),
\end{equation}
where $x_{i}\in \mathbb{R}^{m}$ are the interpolation points, $\psi$ is a radial function and  $w_{i}$ are proper weights. Moreover, the function $p_{j}(x)$ outputs the $j$-th coordinate of $x$. 
The weights $w_i$ are learned by solving the following linear system: 
\begin{equation}
\begin{cases}
\sum_{i=1}^n w_i \psi\left(\left\|x_{k}-x_i\right\|\right)+\sum_{j=1}^{n}z_{j}p_{j}(x_{k})=y_{k} & k=1, \dots, N,\\
\sum_{i=1}^N w_{i} p_{j}(x_{i})=0 & j=1, \dots, m.
\end{cases}
\end{equation}
For small systems, RBFs are faster than neural networks based techniques, but they cubically scale with the number of interpolation points, and may have high memory demands \cite{rbf}.

\subsubsection{Dynamic Mode Decomposition}
Let us suppose to have a sequence of samples $x_{1}, \dots,x_{N_T}$ of a time series with $x_{i}\in \mathbb{R}^{n}.$ 
We want to approximate the Koopman operator $A$: 
\begin{equation}
{x}_{k+1}={A} {x}_k. 
\end{equation}
To achieve this, we build two matrices 
\begin{equation}
S_{1}=[x_{1},\dots,x_{N_T}], \quad
S_{2}=[x_{2},\dots,x_{N_T - 1}].
\end{equation}
We would like to find $A\in \mathbb{R}^{n\times n}$ such that
\begin{equation}
\|{{S}_{2}}-{A} {S}_{1}\|
\end{equation}
is minimized.
However there may be two problems:
\begin{itemize}
    \item the $x_{i}$ may contain some noise, so to solve the system as it is may lead to overfitting problems;
    \item if $n$ is large managing $A$ could be prohibitive.
\end{itemize}
For this reason it is computed a reduced approximation of $A$:
\begin{equation}
\widetilde{A}=UU^{T}S_{2}V\Sigma^{-1}U^{T},
\end{equation}
where $S_{1}=U\Sigma V^{T}$ is a truncated Singular Value Decomposition. There are extensions of DMD that account for the presence of parameter dependent time series\cite{rozza,andreuzzi}.

\subsubsection{Gaussian Process Regression}
A Gaussian process\cite{rozza} is a non parametric model that depends on a positive definite kernel $K.$
Suppose that, given $x_{1},\dots,x_{N}\in \mathbb{R}^{m}$ and $y_{1},\dots,y_{N}\in \mathbb{R}$, we want to model a function $f$ such that $\sum_{i=1}^{n}|f(x_{i})-y_{i}|^{2}$ is minimized.
A Gaussian process assumes that
\begin{equation}
\begin{cases}
y(x)=f(x)+\epsilon, \quad \text{   } \epsilon\sim \mathcal{N}(0,\sigma),\\
[f(x_i)_{i=1 \dots N}] \sim \mathcal{MNN}(g(x_i)_{i=1 \dots N}, K(x_i, x_j)_{i, j=1 \dots N}) ,\\
\end{cases}    
\end{equation}
where $g:\mathbb{R}^{m}\rightarrow \mathbb{R}$ and $\sigma>0.$
It can be proven that for new samples 
$\bar{y}_{i}, i=1\dots l$
we have that
\begin{equation}
\overline{\mathbf{y}} \sim \mathcal{MNN}\left(\hat{\mathbf{K}}\left[\mathbf{K}+\sigma^2 I\right]^{-1} \mathbf{y}, \overline{\mathbf{K}}-\hat{\mathbf{K}}\left[\mathbf{K}+\sigma^2 I\right]^{-1} \hat{\mathbf{K}}^T\right),
\end{equation}
where
\begin{eqnarray}
\overline{\mathbf{y}}=\left\{\bar{y}_i\right\}_{i=1 \dots l}, \\
\mathbf{y}=\left\{y_i\right\}_{i=1 \dots N}, \\
\mathbf{K}=k\left(x_i, x_j\right)_{i, j=1 \dots N}, \\
\hat{\mathbf{K}}=k\left(\bar{x}_i, x_j\right)_{i=1 \dots l, j=1 \dots N}, \\
\overline{\mathbf{K}}=k\left(\bar{x}_i, \bar{x}_j\right)_{i, j=1 \dots l} .
\end{eqnarray}

\subsection{Geometrical Parametrization and Active Subspace Method}
When $\Omega$ depends on $\mu$ the problem becomes more complex because if we want to use the same discretization for every different parameter we need to ensure that it does not produce holes or overlapping elements for every domain. As such, there exists parametrization techniques like Free Form Deformation (FFD)\cite{rozza,siena} or Radial Basis Function (RBF) \cite{rozza,shapeopt3} that are able to parametrize a reference geometry and to morph it without errors.\\
In this context, before proceeding with the development of the surrogate model, it is important to try to reduce the parameter space dimension. At this aim, a technique that can be used is the Active Subspace method\cite{rozza,romor,meneghetti}. 
Let us suppose to have $N$ samples, $x_{1},\dots,x_{N}$, with $x\in \mathbb{R}^{m}$ and a function $f:\mathbb{R}^{m}\rightarrow \mathbb{R}^{n}$ which is also differentiable with continuous derivatives. Then we can construct the matrix $C$ 
\begin{equation}
C=\frac{1}{N}\sum_{i=1}^{N}\nabla f(x_{i})^{T}\nabla f(x_{i}).
\end{equation}
$C$ is symmetric positive definite so it is diagonalizable and admits a decomposition $C=V\Sigma V^{T}$. Assuming that the eigenvalues of $\Sigma$ are sorted in decreasing order one could perform dimensionality reduction by considering $z_{i}=Wx_{i}$ where $W$ are the first $k$ columns of $V$.

\section{Academic benchmarks}
\label{sec:3}
Now we are going to test the performance of the ROM approaches described in the previous section in terms of efficiency and accuracy against some relevant academic test cases. Hereinafter we omit to explicitly report the dependency of the solution on $\mu$ at the aim to make the notation lighter. 

As first test case we consider a classic fluid dynamics problem, the lid driven cavity\cite{prusak}. The FOM is given by the Navier-Stokes equations endowed with proper initial and boundary conditions:
\begin{equation}
\begin{cases}
\mathbf{u}_{{t}}+ (\mathbf{u}\cdot \nabla)\mathbf{u}- \nu {\Delta} \mathbf{u} + {\nabla} p = 0 & \text { in } \Omega \times (0, T), \\
{\nabla} \cdot \mathbf{u}=0 & \text { in }  \Omega \times (0, T), \\ 
\end{cases}
\end{equation}\label{eq:system}
where $\mathbf{u}$ is the velocity vector field, $p$ is the pressure, and $\nu$ is the kinematic viscosity coefficient.

The computational domain $\Omega$ is $[0, 1]^2$. 
We consider the left corner of the square as the origin of the axes. We impose a no slip boundary condition on the lower and lateral
walls. At the upper wall we prescribe a Dirichlet boundary condition $\mathbf{u}(t, x, 1) = (1,0)$. We consider a structured mesh consisting of $105 \times 105$ square cells. We start the simulations from fluid at rest and we set $\Delta t = 0.005$. 
\begin{figure}[h]
\centering
\includegraphics[width=40mm]{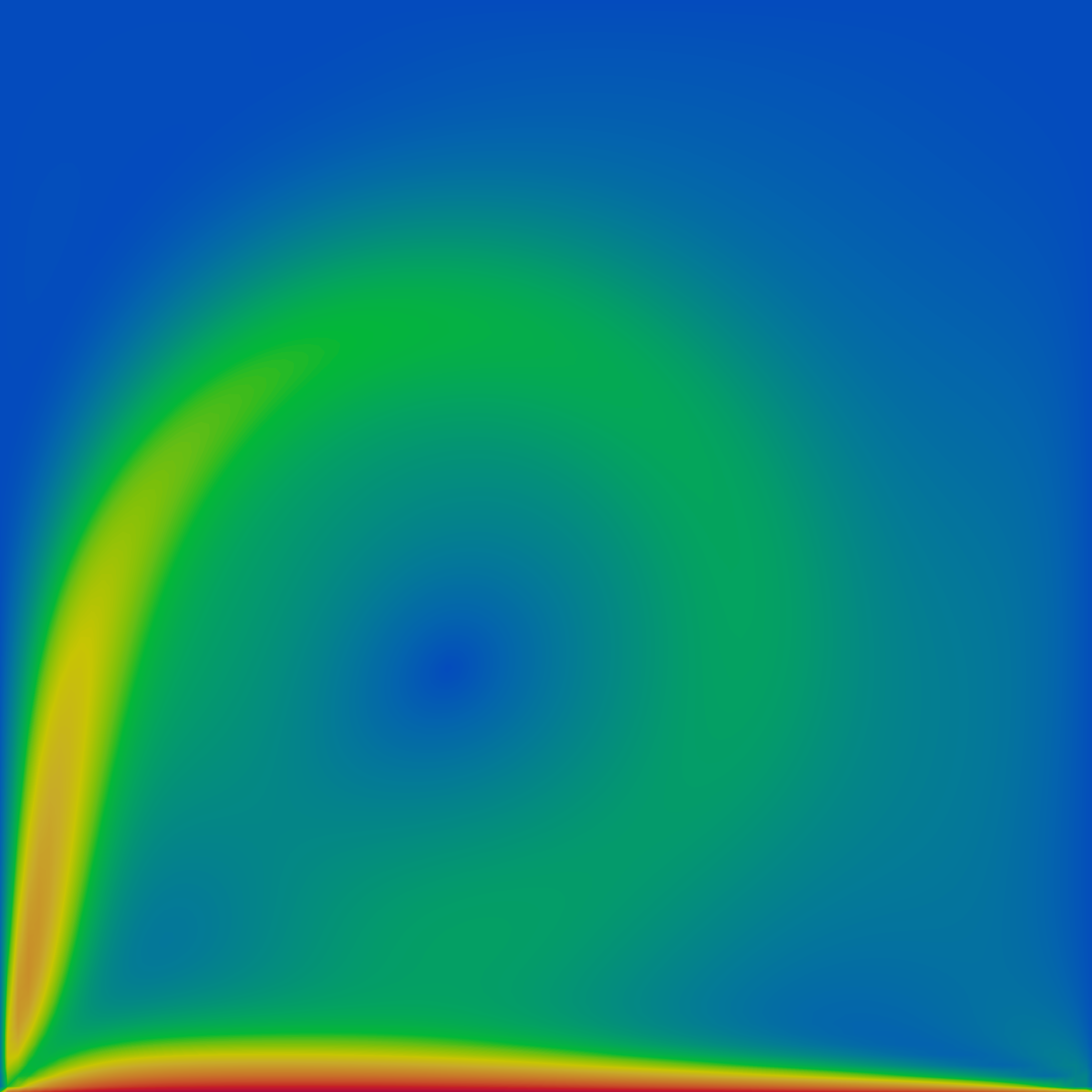}%
\quad  \quad \quad
\includegraphics[width=40mm]{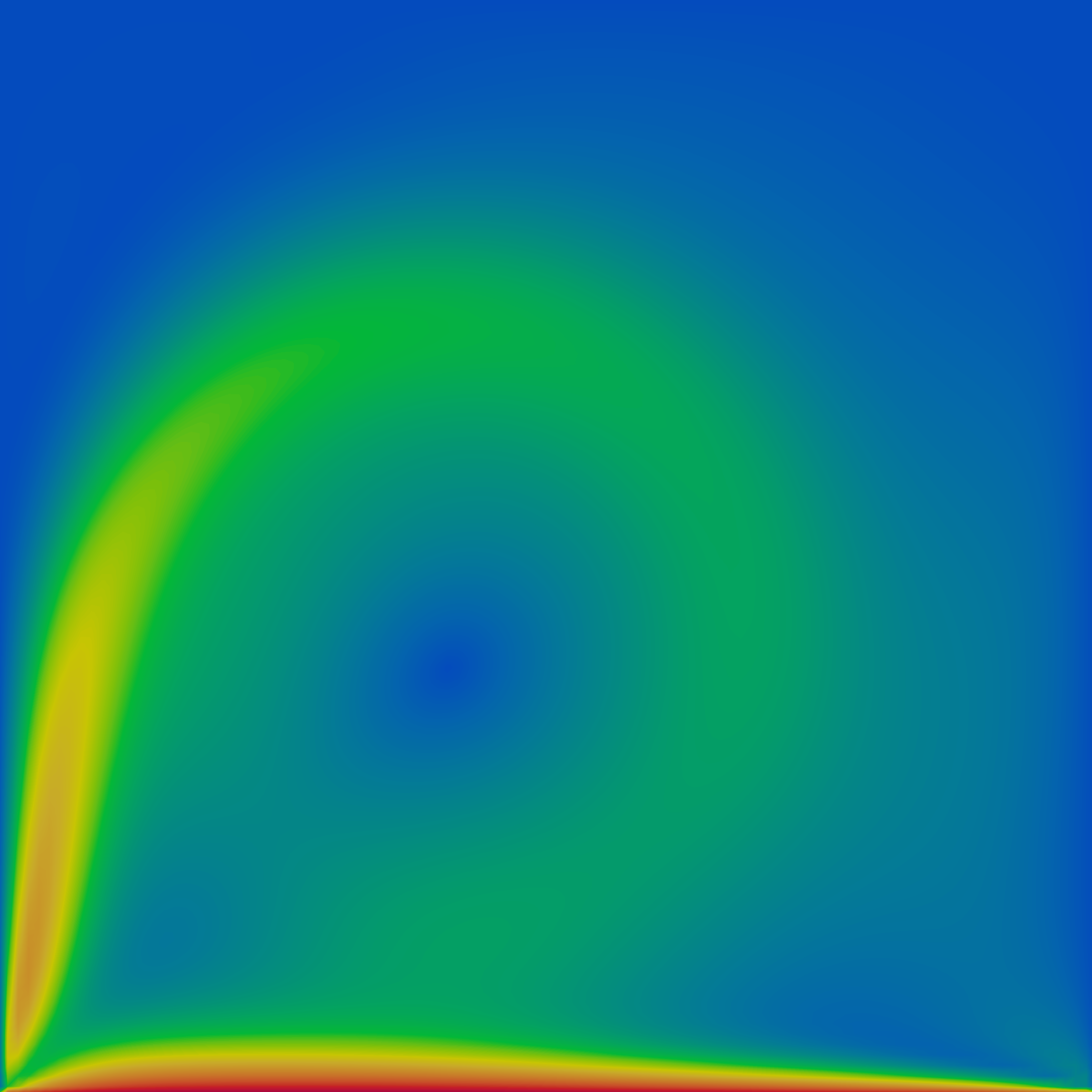}%
\\
\includegraphics[width=65mm]{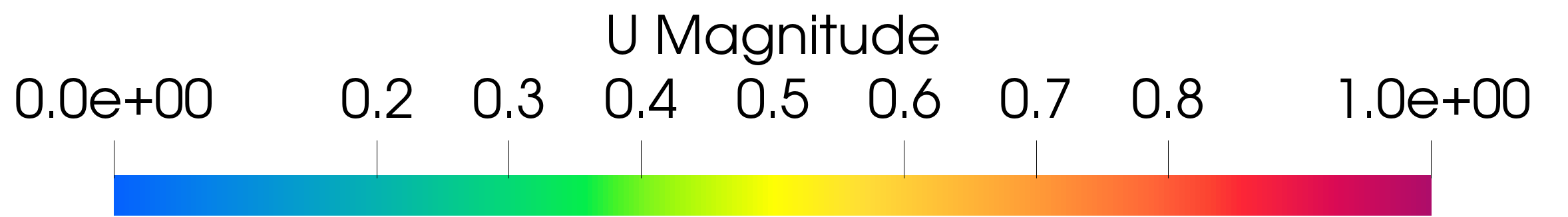}
    \caption{Lid driven cavity: velocity magnitude at $t = 10$ s for $\mu=0.001$ computed by FOM (left) and by DMD+RBF (right).}
\label{fig:lid}
\end{figure}

We set $\mu = \nu$, i.e. the parameter of our problem coincides with the kinematic viscosity. We have  $\mathbb{P} = [0.001,0.01]$.
We consider a uniform sampling of $\mathbb{P}$ consisting of 600 samples to train our ROM. For every training parameter, a FOM simulation was computed by collecting 100 time-dependent snapshots over the interval (0, 10] s. 

All the FOM simulations are performed by using OpenFOAM (\url{https://www.openfoam.com/}), an open source finite volume C++ library widely used by commercial and academic organizations. On the other hand, the ROM
computations are carried out using some python-based packages: PyDMD  (\url{https://pydmd.github.io/PyDMD/}), PyTorch (\url{https://pytorch.org/}) and scikit-learn (\url{https://scikit-learn.org/stable/}).

We take $\mu$ = 0.001 (in the range under consideration but not in the training set) to
evaluate the performance of the parametrized ROM. The results related to the reconstruction of streamwise component of $\mathbf{u}$ are summarized in Table \ref{tab:1}. We consider 400 modes for the DMD. We see that the DMD+RBF model is the most accurate whilst the DDNN one is the most efficient. In Fig. \ref{fig:lid} we show a qualitative comparison between FOM and DMD+RBF.

\begin{table}[h]
\centering
\begin{tabular}{|l|l|l|l|l|}
\hline
                            & DDNN & DMD+DDNN   & DMD+GPR  & DMD+RBF  \\\hline
Train relative error        & 7e-02                     & 2.4e-05                    & 6.9e-05                    & 7.2e-06               \\ \hline
Test relative error         & 6.6e-02                    & 1.8e-05             & 6.3e-05                  &     7.0e-06          \\ \hline
Speedup & 7e4 & 3e2 & 1e2  & 3e2 \\
\hline
\end{tabular}
\caption{Lid driven cavity: performance of ROM models related to the reconstruction of the streamwise component of $\mathbf{u}$. The test value is $\mu = 0.001$}
\label{tab:1}
\end{table}

As a second test case we consider the two dimensional flow past a cylinder\cite{rozza2}. The computational domain is
a 2.2 $\times$ 0.41 rectangular channel with a cylinder of radius 0.05 centered at (0.2, 0.2), when
taking the bottom left corner of the channel as the origin of the axes. 
We impose a no slip boundary condition on the upper and lower wall and on the cylinder. At the inflow and the outflow we prescribe the following velocity
profile 
\begin{equation}
\mathbf{u}(t,0,y) = \left(\dfrac{6}{0.41^2} y (0.41 - y), 0\right), \quad y \in  [0, 2.2]. 
\end{equation}
We consider a grid consisting of 6000 triangular cells. We start the simulations from fluid at rest and we set $\Delta t = 0.01$.

We set $\mu = \nu$ and $\mathbb{P} = [0.01,0.1]$. We consider a uniform sampling of $\mathbb{P}$ consisting of 600 samples to train our ROM. For every training parameter, a FOM simulation was computed by collecting 100 time-dependent snapshots over the interval (0, 8] s. 

All the FOM simulations are performed by using FEniCS (\url{https://fenicsproject.org/}), an open source finite element python library widely used by academic organizations. On the other hand, the ROM computations are carried out using the python packages Rbnics (\url{https://www.rbnicsproject.org/}), PyTorch (\url{https://pytorch.org/}) and scikit-learn (\url{https://scikit-learn.org/stable/}).
\begin{figure}[h]
\centering
\includegraphics[width=100mm]{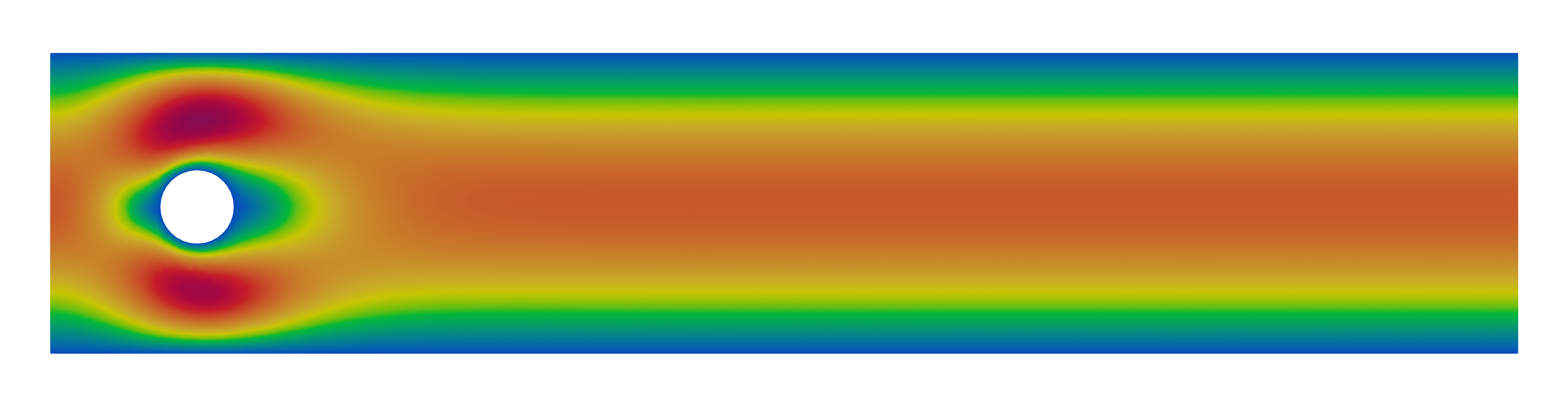}%
\\
\centering
\includegraphics[width=100mm]{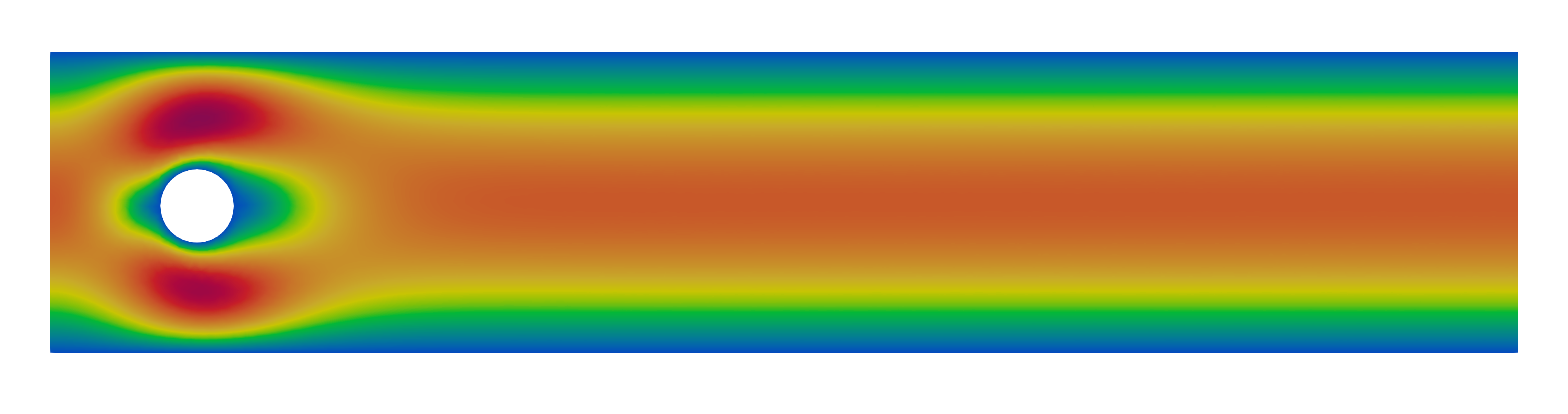}%
\\
\centering
\includegraphics[width=100mm]{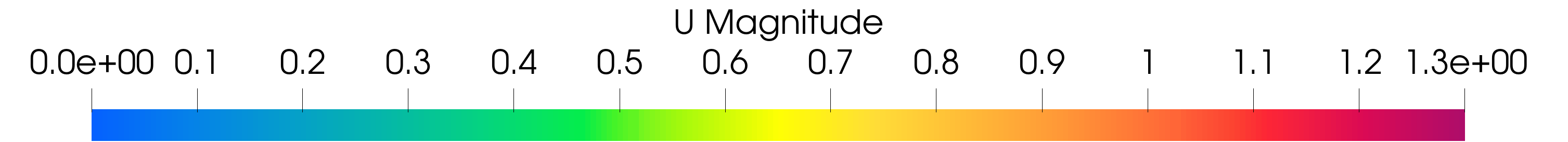}%
\caption{Flow past a cylinder: velocity magnitude at $t = 1$ s for $\mu=0.5$ computed by FOM (top) and by the POD-Galerkin model (with 15 basis functions)}
\label{fig:cyl}
\end{figure}

We take $\mu = 0.5$ (in the range under consideration but not in the training set) to evaluate the performance of the
parametrized ROM. The results related to the reconstruction of the magnitude of $\mathbf{u}$ are summarized in Table \ref{tab:2}. In this case, we see that the Galerkin projection (with 15 basis functions) is the most accurate whilst the RBF one is the most efficient. Moreover it appears evident that, in general, non intrusive ROMs are significantly faster than the Galerkin projection approach although they have in average a lower accuracy. In Fig. \ref{fig:cyl} we show a qualitative comparison between FOM and  Galerkin model for $\mu = 0.5$.

\begin{table}[h]
\centering
\begin{tabular}{|l|l|l|l|l|}
\hline
                            & POD-GALERKIN ($N_{rb}=15$) & POD-GALERKIN ($N_{rb}=5$) & GPR & RBF \\
                            \hline
Train relative error        & 1.6e-06                     & 9.1e-03                    & 1.8e-04                     & 8.8e-17               \\
\hline
Test relative error         & 5.36e-05                    & 3.26e-02                   & 5.05e-03                    & 6.82e-03              \\
\hline
Speedup & 1.5                     & 2          & 3e4    & 1e5    \\
\hline
\end{tabular}
\caption{Flow past a cylinder: performance of ROM models related to the reconstruction of magnitude of $\mathbf{u}$. The test value is $\mu = 0.5$}
\label{tab:2}
\end{table}
As a third and last academic test case we consider an example of geometrical parametrization for a three-dimensional domain. In particular, we refer to the parametrized Stanford Bunny\cite{bunny} problem. We consider the following 
Poisson equation
\begin{equation}
\Delta u = e^{-\left\|x-x_n\right\|^2} \quad \text{in}  \quad \Omega(\mu), \\ 
\end{equation}
where $x_n$ is the vector identifying the barycenter of the domain. We impose the following boundary condition
\begin{equation}
u(x)=e^{-\|x-x_{n}\|} \quad \text{on}  \quad \partial \Omega(\mu).
\end{equation}
The mesh consists of 182184 triangular cells which is parametrized using the Free Form Deformation (FFD) method \cite{rozza}. The parameter $\mu$ of our problem represents the deformation of the control points of the FFD map and we set $\mathbb{P} = [-0.4, 0.4]$. We consider a normal distribution of 600 samples to train our ROM. For every training parameter, the solution of the Poisson equation was computed by using FeniCSx (\url{https://fenicsproject.org/}).

\begin{figure}[h]
\centering
\includegraphics[width=65mm]{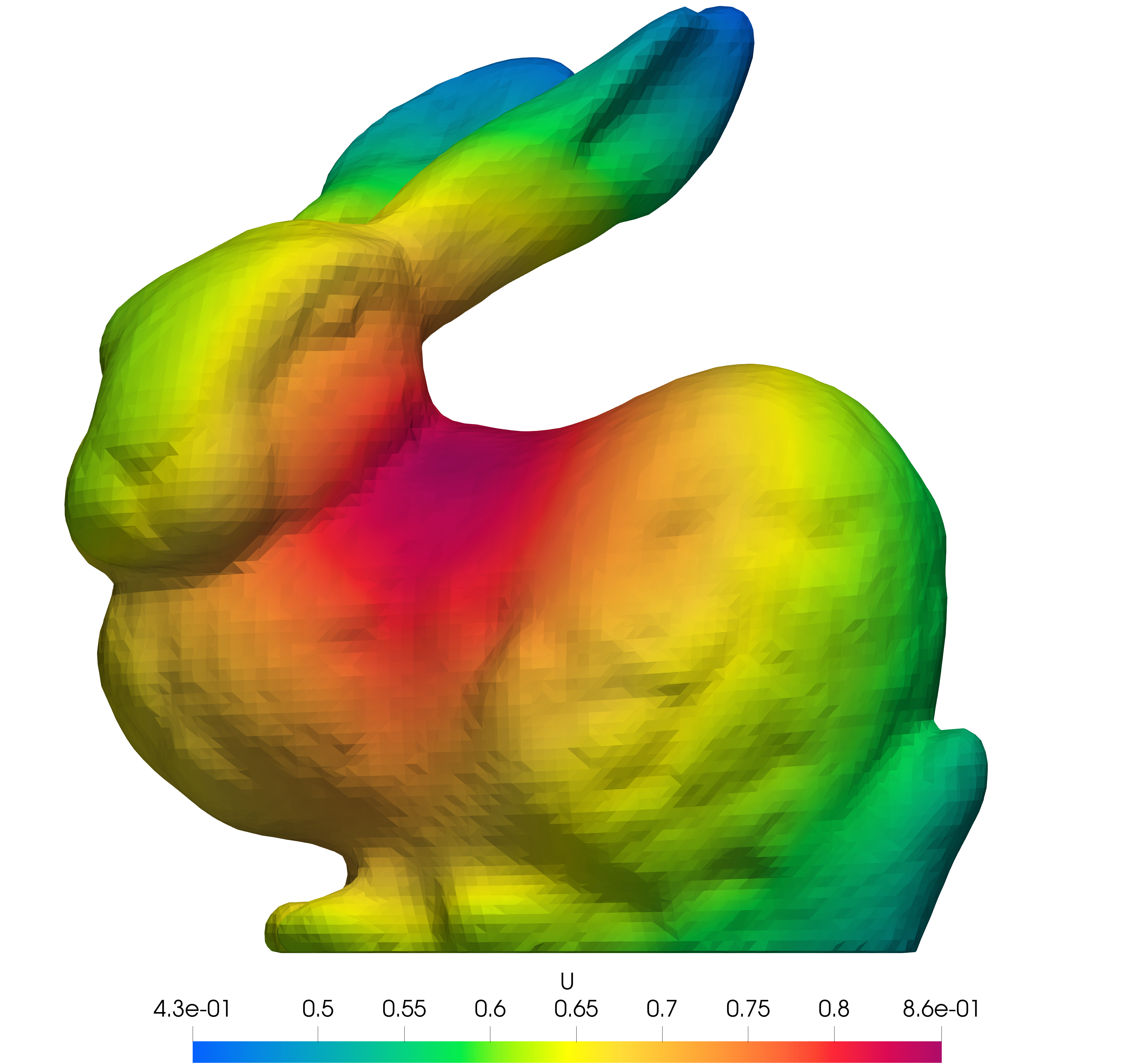}%
\caption{Parametrized Stanford Bunny: plot of the solution for a new value of $\mu$.}
\label{fig:bunny}
\end{figure}

The ROM results are summarized in Table \ref{tab:3}. In this case, we see that the RBF model exhibits the best performance in terms of efficiency and accuracy. In Fig. \ref{fig:bunny} we show the solution obtained for a value of $\mu$ in the range under consideration but not in the training set.

\begin{table}[h]
\centering
\begin{tabular}{|l|l|l|l|}
\hline
                     & RBF    & GPR     & NN   \\
\hline
Train relative error & 4.1e-9 & 7.1e-11 & 4.3e-02 \\
\hline
Test relative error  & 4.7-02 & 5.5e-02 & 4.7e-02 \\
\hline
Speedup & 8e4 & 1.6e4 & 1.6e4 \\
\hline
\end{tabular}\caption{Stanford Bunny test case: performance of ROM models related to the reconstruction of the solution of the Poisson equation on the parametrized domain.}
\label{tab:3}
\end{table}

\section{Towards real world applications: an industrial case}
\label{sec:4}
Now we present an example coming from industry in the framework of a research collaboration with Electrolux concerning domestic refrigerator. Such a problem involves a multiphysics scenario and it is useful to show the performance in terms of accuracy and efficiency of data-driven ROM approaches for complex cases. This case study is still work in progress, so we limit to provide some preliminary outcomes.

At the full order level, air circulation and heat transfer in fluid and between fluid and surrounding solids in a fridge were numerically studied using the Conjugated Heat Transfer (CHT) method to explore both the natural and forced convection-based fridge model followed by a parametric study-based on the ambient temperature $T_{amb}$, fridge fan velocity $\nu_f$, and evaporator temperature $T_{ev}$. The FOM simulations are performed by using OpenFOAM (\url{https://www.openfoam.com}).

A non-intrusive ROM based on a POD-RBF and POD-DDNN methods are considered to obtain the temperature field at specific parametric locations where the training dataset is purely based on numerical computation. The ROM computations were performed in EZyRB (\url{https://github.com/mathLab/EZyRB}). 

To carry out the parametric study, $\nu_f$ ranges from 0\% to 100\% of the full throttle which is 1600
rpm in step of 10 \% i.e., total of 11, two values of $T_{amb}$ i.e., $16^\circ C$ and $32^\circ C$, and four values of $T_{ev}$ i.e., $-15^\circ C$, $-7.9^\circ C$, $-3.25^\circ C$ and $4^\circ C$ are considered.  

\begin{figure}[h]
\centering
\includegraphics[width=150mm]{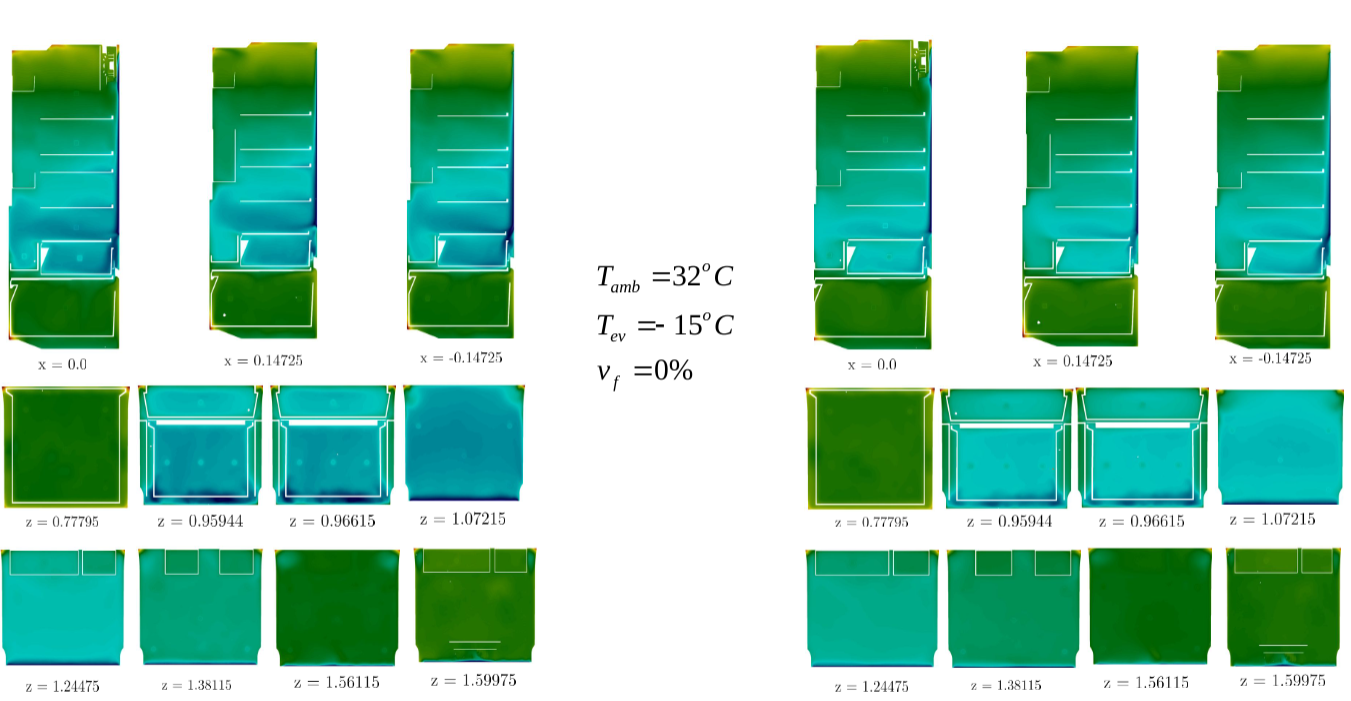}%
\caption{Industrial application: temperature in different sections of the fridge system computed by FOM (left) and by ROM (right). }
\label{fig:fridge}
\end{figure}

In Fig. \ref{fig:fridge} we show the FOM and ROM solutions obtained for $T_{amb} = 32^\circ C$, $T_{amb} = -15^\circ C$ and $\nu_f = 0\%$ in different sections of the fridge system. The comparison is quite satisfactory.  This preliminary result confirms the goodness of data-driven surrogate models resulting to be very promising for complex applications.


\section{Conclusion}
\label{sec:5}
In this paper we have introduced some ROM techniques, both intrusive and non-intrusive. We have tested their performance in terms of efficiency and accuracy against three academic test cases, the lid driven cavity, the flow past a cylinder and the geometrically parametrized Stanford Bunny. Moreover, we also briefly presented some preliminary results related to a more complex case involving an industrial application.

Overall, this contributions shows the applicability of intrusive and non-intrusive ROMs approaches to several contexts.

\section*{Acknowledgements}
We acknowledge the support provided by PRIN “FaReX - Full and Reduced order modelling of coupled systems: focus on non-matching methods and automatic learning” project, INdAM-GNCS 2019–2020 projects, PNRR NGE iN-
EST “Interconnected Nord-Est Innovation Ecosystem” project, and PON “Research and Innovation on Green related issues” FSE REACT-EU 2021 project. This work was also partially supported by the U.S. National Science Foundation through Grant No.
DMS-1953535 (PI A. Quaini).


\begin{thebibliography}{[1]}

\bibitem{democont}
 \textsc{N.~Demo},  \textsc{M.~Strazullo},  and  \textsc{G.~Rozza},
 \jr{Computers \& Mathematics with Applications} \textbf{143}, 383--396 (2023).


\bibitem{digtwin}
 \textsc{M.\,S. Es-haghi},  \textsc{C.~Anitescu},  and
  \textsc{T.~Rabczuk},
 \jr{Computers \& Structures} \textbf{297}, 35 (2024).


\bibitem{padula}
 \textsc{G.~Padula},  \textsc{F.~Romor}, 
 \textsc{G.~Stabile},   and
  \textsc{G.~Rozza},
 \jr{Computer Methods in Applied Mechanics and Engineering} \textbf{423}, 116823 (2024).

\bibitem{papapicco}
 \textsc{D.~Papapicco},  \textsc{N.~Demo}, 
 \textsc{M.~Girfoglio},
 \textsc{G.~Stabile}, 
 and
  \textsc{G.~Rozza},
 \jr{Computer Methods in Applied Mechanics and Engineering} \textbf{393}, 114687 (2022).



\bibitem{shapeopt}
 \textsc{M.~Tezzele}, \textsc{N.~Demo},  \textsc{G.~Stabile},
  \textsc{A.~Mola},  and  \textsc{G.~Rozza},
 \jr{Advanced Modeling and Simulation in Engineering Sciences} \textbf{7}, 19 (2020).

\bibitem{shapeopt3}
 \textsc{A.~Ivagnes},  \textsc{N.~Demo},   and  \textsc{G.~Rozza},
 \jr{International Journal for Numerical Methods in Engineering} \textbf{125}, e7426 (2024).

\bibitem{khamlich}
 \textsc{M.~Khamlich},  \textsc{F.~Pichi},   and  \textsc{G.~Rozza},
 \jr{International Journal for Numerical Methods in Fluids} \textbf{94}, pp.\,1611-1640 (2021).


\bibitem{hajisharifi}
 \textsc{A.~Hajisharifi},  \textsc{M.~Girfoglio},
 \textsc{A.~Quaini},  and  \textsc{G.~Rozza},
 \jr{Finite Elements in Analysis and Design} \textbf{228}, 104050 (2024).


\bibitem{girfoglio}
  \textsc{M.~Girfoglio}, 
 \textsc{A.~Quaini},  and  \textsc{G.~Rozza},
 \jr{Journal of Computational Physics} \textbf{486}, 112127 (2023).

\bibitem{meneghetti}
  \textsc{L.~Meneghetti},  \textsc{N.~Demo},  and  \textsc{G.~Rozza},
 \jr{Applied Intelligence} \textbf{53}, pp.\,22818-22833 (2023).

\bibitem{coscia}
 \textsc{D.~Coscia},  \textsc{N.~Demo},  and  \textsc{G.~Rozza},
 \jr{Scientific Reports} \textbf{14}, 3826 (2024).

\bibitem{balzotti}
 \textsc{C.~Balzotti},  
\textsc{P.~Siena},  
\textsc{M.~Girfoglio},  
    \textsc{G.~Stabile},  
    \textsc{J.~Dueñas-Pamplona},  
    \textsc{J.~Sierra-Pallares},  
 \textsc{I.~Amat-Santos},  and  \textsc{G.~Rozza},
 \jr{ Biomechanics and Modeling in Mechanobiology} (2024).

\bibitem{romor}
 \textsc{F.~Romor},  \textsc{M.~Tezzele},  and  \textsc{G.~Rozza},
 \jr{Journal of Scientific Computing
} \textbf{99}, 83 (2024).

\bibitem{andreuzzi}
 \textsc{F.~Andreuzzi},  \textsc{N.~Demo},  and  \textsc{G.~Rozza},
 \jr{SIAM Journal on Applied Dynamical Systems
} \textbf{22}, pp.\,2432-2458 (2023).

\bibitem{siena}
 \textsc{P.~Siena},
 \textsc{M.~Girfoglio}, \textsc{F.~Ballarin},  and  \textsc{G.~Rozza},
 \jr{Journal of Scientific Computing
} \textbf{94}, 38 (2023).

\bibitem{sheidani}
 \textsc{A.~Sheidani},
 \textsc{S.~Salavatidezfouli},
 \textsc{G.~Stabile},
 \textsc{M.\,B. Gerdroodbary},  and  \textsc{G.~Rozza},
 \jr{Physics of Fluids
} \textbf{35}, 095135 (2023).

\bibitem{pichi}
 \textsc{F.~Pichi},
 \textsc{F.~Ballarin},
 \textsc{G.~Stabile},
 \textsc{G.Rozza},  and  \textsc{J.\,S. Hesthaven},
 \jr{Physics of Fluids
} \textbf{254}, 105813 (2023).

\bibitem{mola}
 \textsc{A.~Mola},
 \textsc{N.~Giuliani},
 \textsc{O.~Crego},
and  \textsc{G.~Rozza},
 \jr{Applied Mathematical Modelling
} \textbf{122}, 322-349 (2023).

\bibitem{prusak}
 \textsc{I.~Prusak},  
 \textsc{D.~Torlo}, 
 \textsc{M.~Nonnino},  and  \textsc{G.~Rozza},
 \jr{Computers \& Mathematics with Applications} \textbf{166}, pp.\,253-268 (2024).

\bibitem{ivagnes}
 \textsc{A.~Ivagnes},  
 \textsc{N.~Tonicello}, 
\textsc{P.~Cinnella}, and
 \textsc{G.~Rozza},
 \jr{arXiv} {2403.05710} (2024).


\bibitem{gonnella}
 \textsc{I.\,C. Gonnella},  
 \textsc{M.\,W Hess}, 
 \textsc{G.~Stabile},  and  \textsc{G.~Rozza},
 \jr{Computers \& Mathematics with Applications} \textbf{149}, pp.\,115-127 (2023).

\bibitem{regazzoni}
 \textsc{F.~Regazzoni},  
 \textsc{S.~Pagani}, 
 \textsc{M.~Salvador},  and  
 \textsc{A.~Quarteroni},
 \jr{Nature Communications} \textbf{15}, 1834 (2024).



\othercit
\bibitem{shapeopt2}
 \textsc{M.~Tezzele},  \textsc{N.~Demo},  and  \textsc{G.~Rozza},
Proceedings of the Viii international conference on computational methods in marine engineering :
  Marine 2019, Gothenburg, Sweden  (International Center for Numerical Methods in Engineering, Barcelona, 2016),  pp.\,122--133.


\othercit
\bibitem{hpc}
 \textsc{F.~Salmoiraghi},  \textsc{F.~Ballarin},  \textsc{G.~Corsi},
  \textsc{A.~Mola},  \textsc{M.~Tezzele},  and  \textsc{G.~Rozza},
Proceedings of the 7th European Conference on Computational Methods in Applied
  Sciences and Engineering, Crete Island, Greece, (Institute of Structural Analysis and Antiseismic Research, Athens, 2016)  pp.\,1013--1031.

\othercit
\bibitem{bunny}
 \textsc{B.~Curless} and  \textsc{M.~Levoy},
Proceedings of the 23rd annual conference on Computer graphics and interactive
  techniques, New Orleans, United States (Association for Computing Machinery, New York, 1996),  pp.\,303--312.


\othercit
\bibitem{hesthaven}
 \textsc{J.~Hesthaven},  \textsc{G.~Rozza},  and  \textsc{B.~Stamm},
Certified Reduced Basis Methods for Parametrized Partial Differential
  Equations, SpringerBriefs in Mathematics (Springer Cham, Berlin, Germany,
  2015).

\othercit
\bibitem{rozza2}
 \textsc{G.~Rozza}, 
 \textsc{F.~Ballarin}, 
  \textsc{L.~Scandurra},  
  and \textsc{F.~Pichi},
Real Time Reduced Order Computational Mechanics,  SISSA Springer Series (Springer Cham, Berlin, Germany,
  2024).


\othercit
\bibitem{hinze}
 \textsc{M.~Hinze},  
 \textsc{J.\,N. Kutz},
 \textsc{O.~Mula},  and  \textsc{K.~Urban},
Model Order Reduction and Applications: Cetraro, Italy 2021, edited by \textsc{M.~Falcone}, and  \textsc{G.~Rozza}, Lecture Notes in Mathematics (Springer Cham, Berlin, Germany,
  2023).

\bibitem{sienalibro}
 \textsc{P.~Siena},  
 \textsc{P.\,C. Africa}, 
\textsc{M.~Girfoglio}, and
 \textsc{G.~Rozza},
 accepted at \jr{Advances in Applied Mechanics} (2024).


\othercit
\bibitem{rozza}
 \textsc{G.~Rozza},  \textsc{G.~Stabile},  and  \textsc{F.~Ballarin},
Advanced Reduced Order Methods and Applications in Computational Fluid Dynamics
  (SIAM, Philadelphia, 2022).

\othercit
\bibitem{rbf}
 \textsc{D.\,H. Buhmann},
Radial Basis Functions: Theory and Implementations (Cambridge University Press,
  Cambridge, 2003).

\othercit
\bibitem{salavatidezfouli}
 \textsc{S.~Salavatidezfouli},
 \textsc{A.~Nikishova},
 \textsc{D.~Torlo},
 \textsc{M.~Teruzzi}, and
 \textsc{G.~Rozza},
Quantitative Sustainability, edited by \textsc{S.~Fantoni}, \textsc{N.~Casagli}, \textsc{C.~Solidoro}, and \textsc{M.~Cobal}  (Springer Cham,
  Berlin, 2024).




\end{thebibliography}
\end{document}